\documentclass[12pt]{article}
\usepackage{amssymb, amsmath}
\usepackage{mathrsfs}       
\usepackage{color}  
     
\numberwithin{equation}{section}

\newtheorem{thm}{Theorem}[section]

\newcommand{\s}{\sigma}

\newcommand{\e}{\varepsilon}
\renewcommand{\a}{\alpha}
\renewcommand{\b}{\beta}
\renewcommand{\d}{\delta}

\renewcommand{\o}{\omega}

\newcommand{\bb}{\begin{equation}}
\newcommand{\ee}{\end{equation}}
\newcommand{\bq}{\begin{eqnarray}}
\newcommand{\eq}{\end{eqnarray}}
\newcommand{\bqn}{\begin{eqnarray*}}
\newcommand{\eqn}{\end{eqnarray*}}
\begin{document}
\title{On a Type I singularity condition in terms of\\  the pressure for the Euler equations in $\mathbb R^3$}
\author{Dongho Chae$^*$  and Peter Constantin$^\dagger$\\
\ \\
 $*$Department of Mathematics\\
Chung-Ang University\\
 Seoul 06974, Republic of Korea\\
 e-mail: dchae@cau.ac.kr\\
and \\
$\dagger$Department of Mathematics\\
Princeton University\\
Princeton, NJ 08544, USA\\
e-mail: const@math.princeton.edu }

\date{}
\maketitle

\begin{abstract}
We prove a blow up criterion  in terms of the Hessian of the pressure of  smooth solutions  $u\in C([0, T); W^{2,q} (\mathbb R^3))$, $q>3$ of the  incompressible Euler equations.   We show that a blow up at $t=T$ happens only if  $$\int_0 ^T \int_0 ^t \left\{\int_0 ^s \|D^2 p  (\tau)\|_{L^\infty}   d\tau \exp \left( \int_{s} ^t  \int_0 ^{\s} \|D^2 p  (\tau)\|_{L^\infty}   d\tau   d\s \right) \right\}dsdt \, = +\infty.$$
As consequences of this criterion we show that there is no blow up  at $t=T$ if  $ \|D^2 p(t)\|_{L^\infty} \le \frac {c}{(T-t)^2}$  with $c<1$ as $t\nearrow T$.  Under the additional assumption of $\int_0 ^T \|u(t)\|_{L^\infty (B(x_0, \rho))} dt <+\infty$, we obtain localized versions of these results.
  \noindent{\bf
keywords:} type I singularity, blow up criterion, Hessian of the pressure 
\end{abstract}

\section{Introduction}
\setcounter{equation}{0}  
We are concerned with the homogeneous incompressible  Euler equation on $\mathbb R^3\times [0, T)$.
$$
(E) \left\{\aligned  &u_t+u\cdot \nabla u= -\nabla p,\label{e1}\\
&\nabla \cdot u=0,\quad u(x,0)=u_0(x)
\endaligned \right.
$$
 where $u(x,t)=(u_1 (x,t), u_2 (x,t), u_3 (x,t))$ is the fluid velocity  and $p=p(x,t)$ is the  scalar pressure. 
 The local in time well-posedness of the Cauchy problem of  (E)  for smooth initial data $u_0$  is well established in various function spaces (see e.g. \cite{maj, con2} and the reference therein). In this paper we are interested in the problem of
  finite time blow up  of local smooth solutions $u\in C([0, T); W^{2, q} (\mathbb R^3))$, for $q>3$, in which cases the local well-posedness is established in \cite{kat}.  There are very many studies of this problem, in particular  establishing blow up criteria (see e.g.  \cite{bkm, koz, con3, con2}).  We mention that there are also  important results  of showing apparition of singularity at  the boundary points  in the domains having boundaries~\cite{luo, elg}.  Our main concern is on the possibility of spontaneous appearance of interior singularity starting from a smooth initial data. 
  Below we first establish a global in space blow up criterion in terms of the Hessian of the pressure, 
  which is sharper than any previously known
  ones in the literature. As an immediate consequence of this criterion we are able to obtain a sharp `small type I condition'  in terms of the Hessian of the pressure, which is consistent with hyperbolic scaling.  The second result is a localization of the first result, 
  showing that certain conditions in terms of the Hessian of the pressure in a ball imply no blow up in the ball.
\begin{thm}
 Let $(u,p)\in C^1 (\mathbb R^3 \times (0, T))$ be a solution of the Euler equation (E) with $ u \in C([0, T); W^{2,q}  (\mathbb R^3))$, for some $ q> 3$.
 \begin{itemize}
\item[(i)] If 
\begin{align}\label{th11}
 \int_0 ^T \int_0 ^t\left\{ \int_0 ^s \|D^2 p  (\tau)\|_{L^\infty}   d\tau \exp \left( \int_{s} ^t  \int_0 ^{\s} \|D^2 p  (\tau)\|_{L^\infty}   d\tau   d\s \right) \right\}dsdt<+\infty, 
\end{align}
then $
 \limsup_{t\to T} \|u(t)\|_{ W^{2,q} } <+\infty.
$
 \item[(ii)]  If 
\bb\label{th12}
 \limsup_{t\to T} \,(T-t)^2 \|D^2 p (t) \|_{L^\infty } <1,
 \ee
 then $
\limsup_{t\to T} \|u(t)\|_{W^{2,q} } <+\infty.
$
 \end{itemize}
  \end{thm}
  \ \\
\noindent{\bf Remark 1.1 } A blow up criterion of the Euler equations in terms of the Hessian of  pressure was obtained  in \cite{cha2} in a different form.  Let $S =(S_{ij})$ with $S_{ij}=\frac12 (\partial_i u_j+\partial_j u_i)$  the symmetric part of the velocity gradient matrix, and we set the unit vectors $\xi= \o/|\o|$,  $\zeta= S\xi /|S\xi|$.  Define $\zeta \cdot P \xi =\mu, $ where 
$P=(\partial_i\partial_j p)$ is the Hessian of the pressure.
 Then, it is shown \cite[Theorem 5.1]{cha2} that there is no blow up at $t=T$ if 
 \bb\label{cpam}
 \int_0 ^T \exp\left(\int_0 ^t \|\mu(s)\|_{L^\infty} ds\right) dt <+\infty.
 \ee
  We note that there exists no mutual implication relation between the conditions \eqref{th11} and \eqref{cpam}, although a stronger condition than \eqref{cpam} (hence, the result of the criterion itself is weaker)
    \bb\label{cpam1}
 \int_0 ^T \exp\left(\int_0 ^t \| D^2 p(s)\|_{L^\infty} ds\right) dt <+\infty.
 \ee  
 implies \eqref{th11}. A related result is found in \cite{con1}.\\
 \ \\
  {\bf Remark 1.2 }  Comparing \eqref{th12} with the `small type I condition'  in terms of the velocity,
  \bb 
   \limsup_{t\to T} \,(T-t) \|Du (t) \|_{L^\infty } <1
   \ee
   introduced in \cite{cha1}, and observing the well-known velocity-pressure relation,    
   $$ \partial_i\partial_j p= \sum_{k,m=1}^3 R_iR_j (\partial_k u_m \partial _m u_k ),
   $$
   we see that \eqref{th12}  is an optimal `small type I condition' in terms of the pressure consistent with hyperbolic scaling. Here, $R_j,$ $j=1,2,3$, are the Riesz tranforms in $\mathbb R^3$(see e.g.\cite{ste})  We also note  that from the condition \eqref{cpam1} it is impossible to deduce the correct small type I condition \eqref{th12} guaranteeing absence of blow up.

\begin{thm}
 Let $(u,p)\in C^1 (B(x_0, \rho) \times (T-\rho, T))$ be a  solution to (E)   with $ u \in C([T-\rho , T); W^{2, q}  ( B(x_0, \rho)))\cap L^\infty ( T-\rho, T; L^2 (B(x_0, \rho)))$ for some $q\in (3, \infty)$. 
  \begin{itemize}

    \item[(i)]  If
 \bb\label{th21}
 \int_{T-\rho} ^{T} \|u(t) \|_{L^\infty (B(x_0, \rho ))} dt<+\infty
 \ee 
and 
 \begin{align}\label{th22} 
& \int_{T-\rho} ^T \left\{\int_{T-\rho} ^t \int_{T-\rho} ^s \|D^2 p  (\tau)\|_{L^\infty (B(x_0, \rho))}   d\tau \times \right.\cr
&\qquad\qquad \left.\exp \left( \int_{s} ^t  \int_{T-\rho} ^{\s} \|D^2 p  (\tau)\|_{L^\infty (B(x_0, \rho))}  d\tau   d\s \right) ds\right\}dt <+\infty,
 \end{align}
then for all $r\in (0, \rho)$
\bb\label{th23}
\limsup_{t\nearrow T} \|u(t)\|_{ W^{2,q} (B(x_0, r))} <+\infty.
\ee
\item[(ii)]  If \eqref{th21} holds,  and 
\bb
\label{th24}
\limsup_{ t\to T} (T-t)^2 \| D^2 p(t) \|_{L^\infty ( B(x_0,\rho))} <1,
\ee
then for all $r\in (0, \rho)$ we have
\bb
\limsup_{t\nearrow T} \|u(t)\|_{ W^{2,q} (B(x_0, r))} <+\infty.
\ee
\end{itemize}
\end{thm}
{\bf Remark 1.3 } A remark similar to Remark 1.2 above holds, comparing \eqref{th24} with the local version of `small Type I condition'
\bb
\limsup_{ t\to T} (T-t) \| \nabla u(t) \|_{L^\infty ( B(x_0,\rho))} <1,
\ee 
 which was obtained in \cite{cw2}.
\section{ The Proof of the Main Theorems}
\setcounter{equation}{0}  
    \noindent{\bf Proof  of Theorem 1.1:  }
     \noindent{\it \underline{Proof of (i)}}  We use the particle trajectory mapping  $\a \mapsto  X ({\a} , t)$ from $\mathbb R^3$ into $\mathbb R^3$   generated by $u=u(x,t)$, which means the solution of  the ordinary differential equation,
\bb\label{particle}
\left\{ \aligned &\frac{\partial   X (\a, t)}{\partial t}   =u( X (\a, t)  , t) \quad \text{on}\quad (0,T),\\
                    & X (\a , 0)  =\a \in \mathbb R^3.\endaligned
                    \right.
\ee
Taking curl of  (E), we obtain  the vorticity equations
   \bb\label{vort}
    \o_t +u\cdot \nabla \o = \o \cdot \nabla u, \quad \o=
    \nabla \times u.\ee
The equation \eqref{vort}  can be rewritten in terms of the particle trajectory as
   \bb\label{21}
    \frac{\partial}{\partial t} \o (X(\a, t), t) = \o (X(\a, t), t)\cdot   \nabla u (X(\a, t), t).
 \ee
    Therefore,
\begin{align}\label{22}
 \frac{\partial^2 }{\partial t^2} \o (X(\a, t), t)   &=   \frac{\partial}{\partial t} \o (X(\a, t), t)   \cdot   \nabla u (X(\a, t), t)\cr
&\qquad+  \o (X(\a, t), t)\cdot   \frac{\partial}{\partial t}    \nabla u (X(\a, t), t).
\end{align}
  From  (E) we also compute
  $$   \frac{\partial}{\partial t}  \partial_j u_k (x,t)  + u\cdot \nabla  \partial_j u_k (x,t)
  =- \sum_{m=1} ^3 \partial_j u_m \partial_m u_k - \partial_j \partial_k p,
  $$
  which can be written as
  \bb\label{23}
  \frac{\partial}{\partial t} \partial_j u_k (X(\a, t), t) 
  =- \sum_{m=1} ^3\partial_j u_m (X(\a, t), t) \partial_m u_k (X(\a, t), t) - \partial_j \partial_k p(X(\a, t), t) .
\ee
  Substituting \eqref{21} and\eqref{23} into \eqref{22}, one has
  \begin{align*}
 \frac{\partial^2 }{\partial t^2} \o_k (X(\a, t), t)   &=  \sum_{j=1}^3 \frac{\partial}{\partial t} \o_j (X(\a, t), t)     \partial_j  u_k (X(\a, t), t)\cr
&\qquad+   \sum_{j=1}^3 \o _j (X(\a, t), t)  \frac{\partial}{\partial t}   \partial_j u_k (X(\a, t), t) 
\end{align*}
  \begin{align}\label{24}
  &=\sum_{j, m=1}^3\o_m (X(\a, t), t)\partial_m u_j (X(\a, t), t)  \partial_j u_k (X(\a, t), t) \cr
&\qquad - \sum_{j, m=1}^3 \o _j (X(\a, t), t)  \partial_j u_m (X(\a, t), t)  \partial_m u_k  (X(\a, t), t)\cr
&\qquad -  \sum_{j=1}^3  \o_j   (X(\a, t), t)  \partial_j\partial_k p  (X(\a, t), t)\cr
& = -  \sum_{j=1}^3  \o_j   (X(\a, t), t)  \partial_j\partial_k p  (X(\a, t), t),
  \end{align}
 from which, after a double integral in time,   we obtain
\begin{align}
\o_k  (X(\a, t), t)&=\o_{0,k} (\a) + \frac{\partial \o_k  (X(\a, t), t)}{\partial t} \Big|_{t=0_+}  t\cr
&\qquad- \sum_{j=1}^3 \int_0 ^t\int_0 ^s  \o_j   (X(\a, \s), \s)  \partial_j\partial_k p  (X(\a, \s), \s)d\s ds\cr
&= \o_{0,k} (\a)+ \sum_{j=1}^3  \o_{0,j } (\a) \partial_j  u_{0, k} (\a) t  \cr
&\qquad - \sum_{j=1}^3 \int_0 ^t\int_0 ^s  \o_j   (X(\a, \s), \s)  \partial_j\partial_k p  (X(\a, \s), \s)d\s ds,
\end{align}
where $\o_0= \nabla \times u_0$, and we used \eqref{21} to compute 
$$\frac{\partial \o_k  (X(\a, t), t)}{\partial t} \Big|_{t=0_+}=  \sum_{j=1}^3  \o_{0,j } (\a) \partial_j  u_{0, k} (\a).
$$  This leads us into
\begin{align}\label{25}
|\o  (X(\a, t), t)| &\le  |\o_0(\a)| + |\o_0(\a) \cdot \nabla u_0 (\a)| t \cr
&\qquad+ \int_0 ^t\int_0 ^s  |D^2 p  (X(\a, \s), \s)|  |\o(X(\a, \s), \s) |d\s ds.
\end{align}
Since the right hand side of \eqref{25} is monotone {non-decreasing} with respect to $t>0$, taking the supremum of the both sides over $(0, t)$, we have
\begin{align}\label{25a}
\sup_{0<\tau<t} |\o  (X(\a, \tau), \tau)|& \le  |\o_0(\a)| + |\o_0(\a) \cdot \nabla u_0 (\a)| t \cr
&\quad+ \int_0 ^t\int_0 ^s  |D^2 p  (X(\a, \s), \s)|  |\o(X(\a, \s), \s) |d\s ds\cr
&\le  |\o_0(\a)| + |\o_0(\a) \cdot \nabla u_0 (\a)| t \cr
&\quad+ \int_0 ^t   \sup_{0<\s<s} |\o  (X(\a, \s), \s)|  \int_0 ^s  |D^2 p  (X(\a, \s), \s)|  d\s ds.\end{align}
Hence the function
$ 
\Phi(s) := \sup_{0<\s <s } |\o(X(\a, \s), \s) |$ satisfies
\begin{align}
 \Phi (t) \le    |\o_0 (\a)|+ |\o_0(\a) \cdot \nabla u_0 (\a)| t+\int_0 ^t  \Phi(s) g(s) ds,  \end{align}
where we set $g(s)=\int_0 ^s |D^2 p  (X(\a, \s), \s)|  d\s$.
 By Gronwall's lemma
 \begin{align}
& \Phi(t)     \le  |\o_0 (\a)| + |\o_0(\a) \cdot \nabla u_0 (\a)| t  \cr
 &\qquad +\int_0 ^t  ( |\o_0 (\a)| + |\o_0(\a) \cdot \nabla u_0 (\a)| s) g(s) \exp \left( \int_{s} ^t  g(\s)   d\s \right) ds\cr
 &   \le \left( |\o_0 (\a)| + |\o_0(\a) \cdot \nabla u_0 (\a)| t \right) \left\{1+ \int_0 ^t   g(s) \exp \left( \int_{s} ^t  g(\s)   d\s \right)ds\right\}.
 \end{align}
 We obtain thus
 \begin{align}\label{keyform}  & |\o(X(\a, t),t)|\le \left( |\o_0 (\a)| + |\o_0(\a) \cdot \nabla u_0 (\a)| t \right)\times\cr
 &\qquad\times \left\{1+\int_0 ^t \int_0 ^s |D^2 p  (X(\a, \tau), \tau)|  d\tau \exp \left( \int_{s} ^t  \int_0 ^{\s} |D^2 p  (X(\a, \tau), \tau)|  d\tau   d\s \right) ds\right\}.\cr
 \end{align}
Taking supremum over $\a \in \mathbb R^3$, and then integrating it over $[0, T]$, we are led to  the inequality
\begin{align}\label{keyforma} 
 & \int_0 ^T \|\o(t)\|_{L^\infty} dt \le \left( \|\o_0 \|_{L^\infty}  + \|\o_0\cdot \nabla u_0 \|_{L^\infty} T \right)\times\cr
 &\qquad\times\left[T+ \int_0 ^T\left\{ \int_0 ^t \int_0 ^s \|D^2 p  (\tau)\|_{L^\infty}   d\tau \exp \left( \int_{s} ^t  \int_0 ^{\s} \|D^2 p  (\tau)\|_{L^\infty}   d\tau   d\s \right) ds\right\}dt\right].\cr
 \end{align}
Applying the well-known Beale-Kato-Majda criterion\cite{bkm} we,  obtain the desired conclusion of (i).\\
 \ \\
   \noindent{\it \underline{Proof of (ii)}}  
   The hypothesis \eqref{th2} implies there exists $ t_0 \in (0, T)$ and $\eta\in (0, 1)$ such that 
 \bb\label{th2}
\sup_{ t_0 < \tau<T } \,(T-\tau)^2 \|D^2 p (\tau) \|_{L^\infty } \le \eta.
\ee 
Hence,
\begin{align}
& \int_{t_0}  ^T \int_{t_0} ^t \left\{\int_{t_0} ^s \|D^2 p  (\tau)\|_{L^\infty}   d\tau \exp \left( \int_{s} ^t  \int_{t_0} ^{\s} \|D^2 p  (\tau)\|_{L^\infty}   d\tau   d\s \right)\right\} dsdt\cr
&\le \eta\int_{t_0}  ^T \int_{t_0} ^t \left\{\int_{t_0} ^s \frac{1}{(T-\tau)^2}   d\tau \exp \left(\eta \int_{s} ^t  \int_{t_0} ^{\s} 
\frac{1}{ (T-\tau)^2}   d\tau   d\s \right) \right\}dsdt\cr
&\le \int_{t_0}  ^T \int_{t_0} ^t\left\{\frac{ 1}{T-s}  \exp \left(\eta \int_{s} ^t  \frac{ 1}{T-\s}d\s \right)\right\} ds dt\cr
&\le\int_{t_0}  ^T\left[ \int_{t_0} ^t\frac{ 1}{T-s}\left(\frac{ T-s}{T-t}\right)^\eta ds\right] dt\cr
&\le \int_{ t_0} ^T \left(\frac{T-t_0}{T-t}\right)^\eta dt= \frac{T-t_0}{1-\eta} <+\infty.
 \end{align}
   The result follows from (i). $\square$\\
   \ \\
   \noindent{\bf Remark:} The starting point of the argument, equation \eqref{24}, can also be derived  from the Lagrangian form of the Euler equations, 
 \bb\label{lagrange} \frac{\partial^2 X(\a, t)}{\partial t^2} =-\nabla p (X(\a, t), t).
 \ee 
Indeed, taking the gradient the both sides of \eqref{lagrange}, and multiplying them by $\o_0 (\a)$, and then using the Cauchy formula
$\o (X(\a, t), t) =\nabla X(\a, t) \o_0 (\a) $, we have \eqref{24}.
   \ \\

   \noindent{ \bf Proof of Theorem 1.2 : }   We consider a sequence $\{t_k\}_{k\in \mathbb N}\subset (T-\rho, T)$ such that $t_k < t_{k+1}$ for all $k\in \mathbb N$,
and $\lim_{k\to +\infty} t_k=T$.  Let  $X(\a, t)$ be the particle trajectory  defined by the ODE in \eqref{particle} for $ (\a, t) \in B(x_0, \rho) \times [T-\rho, T)$.  Thanks to the hypothesis \eqref{th21} we can have a continuous extension of 
   $X(\a, t)$ to an extended domain $B(x_0, \rho) \times [T-r, T]$ by setting $X(\a, T) :=\lim_{t\nearrow T} X(\a, t)$ for all $\a \in B(x_0,\rho)$.
   Indeed, 
  $$ |X(\a, t_k ) -X(\a, t_m)|\le \int_{t_m} ^{t_k}  \|u(t) \|_{L^\infty(B(x_0, \rho)) } dt\to 0
   $$
as $k\ge m\to +\infty$, which shows that for each $\a \in B(x_0, \rho)$ the sequence $\{ X(\a, t_k)\}_{k\in \mathbb N}$ is  Cauchy in $\mathbb R^3$, and converges to a limit.  
Moreover, for we claim    $X(\cdot, T) \in C( B(x_0, r))$ for all $r<\rho$. Indeed, let
$\a, \b\in B(x_0, r)$.  Then,  we have the following estimates for all $\d \in (0, \rho)$.

   \begin{align}\label{th221}  & |X(\a, T) -X(\b, T) |\le |\a-\b| +\int_{T-\rho} ^{T-\d} | u(X(\a, t), t) -u(X(\b, t),t) |dt \cr
  &\qquad+ \int_{T-\e} ^T | u(X(\a, t), t) -u(X(\b, t),t) |dt\cr   
  &\le |\a-\b| + |\a-\b| \int_{T-\rho} ^{T-\d}\frac{ | u(X(\a, t), t) -u(X(\b, t),t) |}{ |X(\a, t)-X(\b, t) |}  \frac{ |X(\a, t) -X(\b, t) |}{|\a-\b|} dt\cr
  &\qquad + 2  \int_{T-\d} ^T( | u(X(\a, t), t)| +|u(X(\b, t),t) | )\cr
      & \le |\a-\b| +|\a-\b| \int_{T-\rho} ^{T-\d} \|\nabla u(t) \|_{L^\infty (B (x_0, r)) } \|\nabla X (\cdot, t)\|_{L^\infty (B (x_0, r)) }  dt \cr
 &\qquad+2 \int_{T-\d} ^T  \|u(t) \|_{L^\infty (B (x_0, r)) }dt\cr
 & \le |\a-\b| \left(1+ \int_{T-\rho} ^{T-\d} \|\nabla u(t) \|_{L^\infty (B (x_0,r)) } e^ {\int_{ T-\rho} ^{t} \|\nabla u(s)\|_{L^\infty (B (x_0, r)) }  ds }dt\right) \cr
 &\qquad+2 \int_{T-\d} ^T  \|u(t) \|_{L^\infty (B (x_0, \rho)) }dt,  
  \end{align}
where we used the estimate
   $$ \| \nabla X(\cdot, t) \|_{L^\infty ( B(x_0, \rho) )} \le  \| \nabla X(\cdot, t_1) \|_{L^\infty ( B(x_0, \rho) )}
   \exp\left( \int^t_{t_1} \| \nabla u(s)\|_{L^\infty ( B(x_0, \rho) )} ds\right),
   $$
which follows from $X(\a, t)=X(\a, t_1) +\int_{t_1} ^ t u (X(\a, s), s) ds$ by taking $\nabla_\a$, and using 
Gronwall's lemma.  We also used the Sobolev inequality
$$ \int_{0} ^{T-\d} \|\nabla u(s) \|_{L^\infty(B(x_0, r)) }  dt\le C(T-\d)\sup_{0<t<T-\d}  \|u(t) \|_{W^{2, q} (B(x_0, \rho)) } <+\infty,
$$
for all $\d \in (0, T)$, which holds thanks to the assumption $u \in C([0, T); W^{2,q}  (B(x_0, \rho)))$ with $ q> 3$.
Now, given $\eta >0 $, we choose $\delta >0$ so that 
$$2 \int_{T-\d} ^T  \|u(t) \|_{L^\infty (B (x_0, \rho)) }dt <\frac{\eta}{3},$$
Then, for such  $\delta >0$  we choose $|\a-\b|$ small enough  to have
$$
|\a-\b| \left(1+ \int_{0} ^{T-\d} \|\nabla u(t) \|_{L^\infty (B (x_0,r)) } e^ {\int_{ T-\rho} ^{t} \|\nabla u(s)\|_{L^\infty (B (x_0, r)) }  ds }dt \right)<\frac{\eta}{3}.
$$
Then, \eqref{th221} shows that $ |X(\a, T) -X(\b, T) |< \eta$. The claim  $X(\cdot, T) \in C( B(x_0, r))$ is proved.\\

By the continuity of the trajectory mapping $ X(\cdot , t)$ for $t\in [0, T]$  for each $r\in (0, \rho)$ there exists $\e >0$ such that 
 \bb
B\left(x_0, r+\frac{\rho-r}{3} \right)\subset  X\left( B\left(x_0, r+\frac{\rho-r}{2} \right), T-t\right) \subset  B(x_0, \rho)\quad \forall t\in [0, \e].
\ee
   where $ X(\a, t) $ is the extension of  the particle trajectory  to $t=T$ defined  by the following ordinary differential equations
 \bb\label{particle1}
\left\{ \aligned &\frac{\partial   X (\a, t)}{\partial t}   =u( X (\a, t)  , t) \quad \text{on}\quad [T-\e, T) ,\\
                    & X (\a , T-\e)  =\a \in  B(x_0, r).
                    \endaligned
                    \right.
\ee
   Then,
   we have from \eqref{keyform} 
  \begin{align}
 & \|\o (t) \| _{ L^\infty (B\left(x_0, r+\frac{\rho-r}{3} \right)) } \le \sup_{\a\in B\left(x_0, r+\frac{\rho-r}{2} \right)}  \left\{  | \o (\a,  T-\e ) |  + | \o (\a, T-\e)||\nabla u(\a,  T-\e ) t \right\}  \times \cr
     &\quad\times\left[1+\int_{T-\e}^t \left\{  \int_{T-\e} ^s |D^2 p  (X(\a, \tau), \tau)|  d\tau \exp \left( \int_{s} ^t  \int_{T-\e}^{\s} |D^2 p  (X(\a, \tau), \tau)|  d\tau   d\s \right)ds\right\}\right]\cr
   & \le  \left\{  \| \o_0\|_{L^\infty (B(x_0,\rho))}   + \| \o_0 \nabla u_0\|_{L^\infty (B(x_0,\rho))}  T\right\}  \times \cr
     &\quad\times\left[ 1+\int_{T-\e}^t \left\{\int_{T-\e} ^s \|D^2 p  (\tau)\|_{L^\infty (B(x_0,\rho))}   d\tau \exp \left( \int_{s} ^t  \int_{T-\e}^{\s} \|D^2 p  ( \tau)\|_{L^\infty (B(x_0,\rho))}   d\tau   d\s \right)ds\right\} \right]. 
     \end{align}
Integrating this over $[T-\e, T]$, we find
   \begin{align*}
 & \int_{T-\e} ^T \|\o (t) \| _{ L^\infty (B\left(x_0, r+\frac{\rho-r}{3} \right)) }   dt  \le  \left\{  \| \o_0\|_{L^\infty (B(x_0,\rho))}   + \| \o_0 \nabla u_0\|_{L^\infty (B(x_0,\rho))}  T\right\}  \times \cr
     &\quad\times\left[T+\int_{T-\e} ^T\int_{T-\e}^t \left\{ \int_{T-\e} ^s \|D^2 p  (\tau)\|_{L^\infty (B(x_0,\rho))}   d\tau\times \right. \right.\cr
     &\hspace{2.in}\quad \left.\left. \times \exp \left( \int_{s} ^t  \int_{T-\e}^{\s} \|D^2 p  ( \tau)\|_{L^\infty (B(x_0,\rho))}   d\tau   d\s \right) ds\right\} dt \right],
      \end{align*}
  which is finite by the hypothesis \eqref{th22}.  Applying the localized version Beale-Kato-Majda criterion  of \cite[Theorem 1.1]{cw1}, we obtain 
   \eqref{th23}. This completes the proof of (i).\\
The proof of (ii) is  exactly the same as  that of Theorem 1.1(ii) above, just replacing $ \|D^2 p(t) \|_{L^\infty} $ by $\| D^2 p(t) \|_{L^\infty ( B(x_0,\rho))}$. $\square$\\
\ \\
$$ \text{\bf Acknowledgements} $$
Chae was partially supported by NRF grant 2016R1A2B3011647, while the work of P.C. was partially supported by the Simons Center for Hidden Symmetries and Fusion Energy. The authors declare  that there is no conflict of interest.

\end{document}